\documentclass[11pt]{article}

\usepackage{amsfonts,amsmath,latexsym,color,epsfig,enumitem, amssymb,bbm,amsthm}
\usepackage{hyperref}
\usepackage{subcaption}
\usepackage{ulem}

\newtheorem{theorem}{Theorem}[section]

\newtheorem{lemma}{Lemma}[section]

\theoremstyle{definition}

\numberwithin{equation}{section}

\renewcommand{\le}{\leqslant}
\renewcommand{\ge}{\geqslant}
\renewcommand{\leq}{\leqslant}
\renewcommand{\geq}{\geqslant}
\newcommand{\R}{\mathbb R}

\newcommand{\E}{\mathbb E}

\def\qed{\ifvmode\mbox{ }\else\unskip\fi\hskip 1em plus 10fill$\Box$}

\input{epsf}

\makeatletter
\def\Ddots{\mathinner{\mkern1mu\raise\p@
\vbox{\kern7\p@\hbox{.}}\mkern2mu
\raise4\p@\hbox{.}\mkern2mu\raise7\p@\hbox{.}\mkern1mu}}
\makeatother

\def\R{\mathbb R}

\def\E{\mathbb E}

\title{\vspace{-0.7cm}On the number of high-dimensional partitions}
\author{Cosmin Pohoata\thanks{Department of Mathematics, Emory University, Atlanta, GA. Email: {\tt cosmin.pohoata@emory.edu}. Research supported by NSF Award DMS-2246659.} \and Dmitriy Zakharov\thanks{Department of Mathematics, Massachusetts Institute of Technology, Cambridge, MA. Email: {\tt  zakhdm@mit.edu}. Research supported by the Jane Street Graduate Fellowship.}}
\date{}

\begin{document}
\maketitle

\begin{abstract}
Let $P_{d}(n)$ denote the number of $n \times \ldots \times n$ $d$-dimensional partitions with entries from $\left\{0,1,\ldots,n\right\}$. Building upon the works of Balogh-Treglown-Wagner and Noel-Scott-Sudakov, we show that when $d \to \infty$,
$$P_{d}(n) = 2^{(1+o_{d}(1)) \sqrt{\frac{6}{(d+1)\pi}} \cdot n^{d}}$$
holds for all $n \geq 1$. This makes progress towards a conjecture of Moshkovitz-Shapira [{\it{Adv. in Math.}} 262 (2014), 1107--1129]. Via the main result of Moshkovitz and Shapira, our estimate also determines asymptotically a Ramsey theoretic parameter related to Erd\H{o}s-Szekeres-type functions, thus solving a problem of Fox, Pach, Sudakov, and Suk [{\it{Proc. Lond. Math. Soc.}} 105 (2012), 953--982]. Our main result is a new supersaturation theorem for antichains in $[n]^{d}$, which may be of independent interest. 


\end{abstract}

\section{Introduction}

Counting partitions is arguably one of the most important topics in the development of modern combinatorics and number theory. Nevertheless, while thousands of papers have been written about it by now, many interesting questions still remain open. A few particularly beautiful ones deal with so-called high-dimensional partitions, which still seem to continue to be rather mysterious objects after all these years. A $d$-dimensional partition of a positive integer $N$ is typically defined as a $N \times \ldots \times N$ $d$-dimensional tensor $A$ with nonnegative integer entries which sum up to $N$, and so that
$$A_{i_1,\ldots,i_{t},\ldots,i_{d}} \geq A_{i_1,\ldots,i_{t}+1,\ldots,i_{d}}$$
holds for all possible indices $i_{1},\ldots,i_{d}$ and for all $t \in \left\{1,\ldots,d\right\}$. When $d = 1$, this is simply a decreasing sequence of nonnegative integers $A_{1} \geq A_{2} \geq \ldots$, also known sometimes in this context as a {\it{line partition}}. Entries equal to zero are traditionally omitted. In particular, if $p_{d}(N)$ denotes the total number of $d$-dimensional partitions of $N$, then $p_{1}(N)$ is simply the usual partition number $p(N)$. 

An intimately related parameter is the number $P_{d}(n)$ of $n \times \ldots \times n$ $d$-dimensional partitions with entries from $\left\{0,1,\ldots,n\right\}$. These are also $d$-dimensional tensors $A$ with nonnegative integer entries, so that $A_{i_1,\ldots,i_{t},\ldots,i_{d}} \geq A_{i_1,\ldots,i_{t}+1,\ldots,i_{d}}$
holds for all possible indices $i_{1},\ldots,i_{d}$ and for all $t \in \left\{1,\ldots,d\right\}$, like before, but without a fixed sum of entries. Instead, the range for the entries is the interval $\left\{0,1,\ldots,n\right\}$. When $d = 1$, this is typically visualized as a $2$-dimensional sequence of stacks of height $A_i$ each, which is commonly known as the Young diagram of the line partition. 

With this perspective, it is a simple counting exercise to see that
$$P_{1}(n) = {2n \choose n}.$$
However, it turns out that computing $P_{d}(n)$ when $d \geq 2$ is substantially more difficult. A celebrated result of MacMahon \cite{M60} states that
$$P_{2}(n) = \prod_{1 \leq i,j,k \leq n} \frac{i+j+k-1}{i+j+k-2},$$
but finding a similar formula when $d \geq 3$ is an outstanding open problem in enumerative combinatorics. In fact, for $d \geq 3$, even the generating functions of $P_{d}(n)$ and $p_{d}(N)$ are already unknown, the only available results being the classical ones by Euler and MacMahon for $d=1$ and $d=2$, respectively. We refer to \cite{Andrews} and \cite{Kratt} for more background and references on high-dimensional partitions. 

In this paper, we settle a problem raised by Moshkovitz and Shapira in \cite{MShap}, regarding the asymptotic behavior of $P_{d}(n)$ when $d \to \infty$. This turns out to be much more difficult to understand than that of $p_{d}(n)$ in this regime (one can easily prove that $p_{d}(n) \sim {d \choose n-1}$ for $d\rightarrow \infty$ and fixed $n$).
Our main result is the following estimate.

\begin{theorem} \label{main}
For every $n \geq 1$,
$$P_{d}(n) = 2^{(1+o_{d \to \infty}(1)) \sqrt{\frac{6}{(d+1)\pi}} \cdot n^{d}}.$$
\end{theorem}

We deduce Theorem \ref{main} from a more precise result about the number of antichains in the $(d+1)$-dimensional box $[n]^{d+1}$, where $[n]=\left\{1\ldots,n\right\}$. To set it up, we first briefly recall the natural poset structure on $[n]^{d}$ and some important standard results about it.

For $x,y \in [n]^{d+1}$, we write $x \preccurlyeq y$ whenever $x_{i} \leq y_{i}$ holds for all $1 \leq i \leq d+1$. Whenever we know that at least one of these $d+1$ inequalities is strict, we shall simply write $x \prec y$.  A subset $A \subset [n]^{d+1}$ is an {\it{antichain}} if no two elements $x \neq y$ satisfy $x \preccurlyeq y$ or $y \preccurlyeq x$. Or in other words, $A$ is an antichain if for every $a \in A$, $x \not\in A$ for every element $x \in [n]^{d+1}$ such that $x \preccurlyeq a$. It is a well-known theorem of de Bruijn et al. \cite{BTK51} that the poset $[n]^{d+1}$ equipped with $\preccurlyeq$ satisfies the {\it{Sperner property}}, which means here that the largest antichain in $[n]^{d+1}$ is also the so-called ``middle layer". In this case, this is the set of solutions $(x_{1},\ldots,x_{d+1}) \in [n]^{d+1}$ to the equation $x_{1}+\ldots+x_{d+1} = (d+1)(n+1)/2$, whose cardinality we denote by $N_{d+1,n}$. It is worth noting that for $n=2$ this is precisely the statement of the classical Sperner's theorem \cite{Sperner} in disguise. See, for example, \cite{KG} for a bit more about the history of these results. For arbitrary $n$, using a version of the central limit theorem, it can be shown that, when $d \to \infty$, we have that
\begin{equation} \label{CLT} 
N_{d+1,n} = \sqrt{ \frac{6}{\pi(n^2-1)(d+1)}} \cdot n^{d+1} ( 1+o_{d \to \infty}(1)) \sim \sqrt{ \frac{6}{(d+1) \pi} } \cdot n^{d}.
\end{equation}
For more details, see \cite{MR08}. We are now ready to state our more technical result.

\smallskip

\begin{theorem} \label{anti}
For every $n \geq 2$, the total number $\mathcal{A}_{d+1,n}$ of antichains in the $(d+1)$-dimensional box $[n]^{d+1}$ satisfies
$$\mathcal{A}_{d+1,n} \leq 2^{(1 + \frac{C_n \log d}{\sqrt{d}}) N_{d+1,n}},
$$
where $C_n > 0$ is a constant which solely depends on $n$.
\end{theorem}

\smallskip

A set $B \subset [n]^{d+1}$ is called a {\it{down-set}} if $b \in B$ implies $x \in B$ for all $x \preccurlyeq b$. The number of such objects is also related to the number of $d$-dimensional partitions with entries in $\left\{0,1,\ldots,n\right\}$. In fact, the connection between these different worlds is rather simple, as also noted, for example, in \cite{MShap}. Every down-set $B$ comes with a natural $d$-dimensional partition with entries in $\left\{0,1,\ldots,n\right\}$ defined as follows: for every $1 \leq i_1,\ldots,i_{d} \leq n$, let 
$$A_{i_1,\ldots,i_{d}} = \max\{ s\ :\ (i_1,\ldots,i_d,s) \in B),$$
whereas if this set is ever empty, simply set $A_{i_1,\ldots,i_{d}}=0$. Since $B$ is a down-set, the condition $A_{i_{1},\ldots,i_{t},\ldots,i_{d}} \geq A_{i_{1},\ldots,i_{t}+1,\ldots,i_{d}}$ is satisfied or every $1 \leq t \leq d$. Furthermore, it is easy to check that this map is a bijection. At the same time, down-sets in $[n]^{d+1}$ are also in a natural one-to-one correspondence with antichains in $[n]^{d+1}$: the $\preccurlyeq$-maximal elements of $B$ yield an antichain, and it is easy to see that every antichain can be obtained from a down-set in this way. 

This then implies that $\mathcal{A}_{d+1,n}=P_{d}(n)$, and so, by Theorem \ref{anti} and the estimate from \eqref{CLT}, we get that
$$\log_{2} \mathcal{P}_{d}(n) \leq \left(1 + \frac{C_n \log d}{\sqrt{d}}\right) N_{d+1,n} = (1+o_{d \to \infty}(1)) \sqrt{\frac{6}{(d+1)\pi}} \cdot n^{d}.$$
 
In order to prove Theorem \ref{main}, it therefore suffices to prove the estimate from Theorem \ref{anti}. 

We will start discussing our proof of Theorem \ref{anti} in the next section, but for now it is important to at least mention that Theorem \ref{anti} rests upon a few essential previous works. First of all, when $n=2$, Theorem \ref{anti} recovers a celebrated result of Kleitman from 1969, which at the time settled a longstanding problem of Dedekind. See \cite{Kleitman} and the many references therein. Second, when $n=3$, Theorem \ref{anti} also recovers a recent difficult extension of Kleitman's result by Noel, Scott and Sudakov \cite{NSS}. The latter builds upon another remarkable recent paper of Balogh, Treglown, and Wagner \cite{BTW}, who provided an alternative proof of Kleitman's theorem using the so-called graph container method, originally introduced by Kleitman and Winston in \cite{KW1} and \cite{KW2} (and at the heart of several other exciting developments in extremal combinatorics in the last few years, see \cite{BMS} and \cite{ST}). The proof of Theorem \ref{anti} also relies on this alternative approach, our main input essentially consisting of a new supersaturation result for antichains in $[n]^{d+1}$, which we establish in full generality for all $n \geq 2$.

Last but not least, it is perhaps also important to mention that via results of Moshkovitz and Shapira from \cite{MShap}, Theorem \ref{anti} also completely determines, in the $d \to \infty$ limit, an elegant Ramsey theoretic parameter $N_{3}(d+1,n)$ related to Erd\H{o}s-Szekeres-type functions, originally introduced by Fox, Pach, Suk, and Sudakov in \cite{FPSS}. This is formally defined as the smallest integer $N$ so that every $(d+1)$-coloring of the edges of the complete $3$-uniform hypergraph on $\left\{1,\ldots,N\right\}$ must always yield a monochromatic monotone path of length $n$. By \cite[Theorem 1]{MShap}, it turns out that $N_{3}(d+1,n) = P_{d}(n)+1$ holds, so the estimate from Theorem \ref{main} completely solves the first open problem raised in the last section of \cite{FPSS}. 

\section{Proof of Theorem \ref{anti}}

For convenience, let us denote $d+1$ by $D$, and translate the box to include the origin, which means that from now on $[n]$ will stand for $\left\{0,1,\ldots,n-1\right\}$ instead of $\left\{1,\ldots,n\right\}$. The poset structure on $[n]^{D}$ is of course defined the same way: for $x,y \in [n]^{D}$, recall that we have $x \preccurlyeq y$ whenever $x_{i} \leq y_{i}$ holds for all $1 \leq i \leq D$. However, it is worth emphasizing that the width of the middle layer of $[n]^{D}$ under this transformation gets renormalized. We shall thus now denote by $N_{n,D}$ the number of sequences $(a_1, \ldots, a_D)$ such that $a_i \in [n]$ and $\sum_{i=1}^{D} a_i = \frac{n-1}{2}D$. With this notation, the Sperner property of $[n]^{D}$ becomes the statement that the largest antichain in $[n]^{D}$ is of size $N_{n,D}$. 

For every $n \geq 2$, it therefore suffices to show that the total number $\mathcal{A}_{D,n}$ of antichains in the $D$-dimensional box $[n]^{D}$ satisfies
\begin{equation} \label{anti2}
\mathcal{A}_{D,n} \leq 2^{(1 + \frac{C_n \log D}{\sqrt{D}}) N_{n,D}},
\end{equation}
where $C_n > 0$ is a constant which solely depends on $n$. 

We prove this using the so-called graph container method, which we will apply in the style of Noel, Scott, and Sudakov \cite{NSS}. Our main new ingredient is the following supersaturation result.

\begin{theorem}\label{super}
Let $a$ and $b$ be positive integers and let $X \subset [n]^D$ be a set of size at least $a (1+ 3n/D) N_{n, D} + b$ for some $a$ and $b$. Then $X$ contains at least $\frac{b}{a! n^a}  D^a$ comparable pairs.
\end{theorem}

We first discuss the proof of Theorem \ref{super} in Section 2.1 below, and then use this result to derive the estimate from \eqref{anti2} in Section 2.2. 

\subsection{Supersaturation in $[n]^{D}$} 

We shall consider the natural grading of $[n]^{D}$ induced by the rank function $r : [n]^D \to \mathbb{N}$ defined by 
$$r(x) = x_1+\ldots+x_D,\ \ \text{for every}\ x = (x_1,\ldots,x_D) \in [n]^D.$$
This gives rise to the partition 
$$[n]^{D} =\bigsqcup_{i=0}^{(n-1)D} V_i,$$
into ``level sets" $V_{i} = \left\{x \in [n]^D:\ r(x) = i\right\}$, and it is easy to see that $V_1,\ldots,V_{(n-1)D}$ indeed satisfy the following two nice properties:

\begin{itemize}
    \item if $x \in V_i$ and $x \prec y$, then $y \in \bigcup_{j > i} V_j$,
    \item if $x \prec y$ are such that $x \in V_i$ and $y \in V_j$ for $j > i+1$, then there exists $z \in V_{i+1}$ such that $x<z<y$.
\end{itemize}

Our strategy to prove Theorem \ref{super} will be roughly as follows. For each $i = 0, \ldots, (n-1)D-1$, we consider a weighted version of the comparability graph between $V_i$ and $V_{i+1}$, and choose a random (partial) matching $M_i$ so that each edge in this graph is in $M_i$ with a prescribed probability (depending on its weight). To do this, we construct a suitable probability distribution $\mu_i$ on the set of matchings between $V_i$ and $V_{i+1}$ and we shall sample $M_i$ according to $\mu_i$. In the spirit of Stanley \cite{StanleyHL}, we shall then consider the set of chains $\{C_1, \ldots, C_N\}$ determined by the union of these matchings. The chains are such that $[n]^D= \bigsqcup_{j = 1}^N C_j$, and, as we will show, their number is not much larger than the size of the middle layer. More precisely, $N \le (1 + 3n/D) N_{n, D}$ will hold true. Consequently, by the pigeonhole principle, every set $X \subset [n]^D$ of size at least $a(1 + 3n/D) N_{n, D} + b$ will always come with at least $b$ pairs $(x, y)$ of elements of $X$ belonging to the same set $C_j$ and such that $x \preccurlyeq y$, and the rank of $y$ exceeds the rank of $x$ by at least $a$. To get a reasonable lower bound for the number of comparable pairs in $X$, one can then hope to proceed by estimating, for each fixed pair $x \preccurlyeq y$ with $x, y \in X$, the probability that $x, y$ belong to the same set $C_j$ for some $j$.

\paragraph{Distributions on matchings.} To construct our desired probability distributions, we shall rely on the following general observation.
Let $G= (A \cup B, E)$ be a bipartite multigraph with parts $A$ and $B$. Let $\mathcal M_m$ be the family of all matchings between $A$ and $B$ in $G$ of size $m$. 
\begin{lemma} \label{polytope}
Suppose that $f: E \rightarrow [0, 1]$ is such that for any vertex $v \in A$ we have $\sum_{u \in B} f(v,u) \le 1$ and for any $u \in B$ we have $\sum_{v \in A} f(v,u) \le 1$. Suppose that $\sum_{(v,u) \in E} f(v,u) = m$ is an integer.
Then, there always exists a probability distribution $\mu: \mathcal M_m \rightarrow [0, 1]$ such that a $\mu$-random matching $M \in \mathcal M_m$ contains a fixed edge $e \in E$ with probability $f(e)$. 
\end{lemma}

From now on, we will typically use $\sum_{v \sim u}$ instead of $\sum_{(v,u) \in E}$ for every summation over all the edges in $E$, whenever the underlying graph is clear from the context. 

\begin{proof}
Let $P \subset [0, 1]^E \subset \R^E$ be the set of all functions $g: E \rightarrow [0, 1]$ such that for any $v \in A$ we have $\sum_{u \in B} g(v, u) \le 1$, for any $u \in B$ we have $\sum_{v \in A} g(v, u) \le 1$ and $\sum_{v \sim u} g(v, u) = m$. If we fix a linear order $e_1,\ldots,e_{|E|}$ on the set of edges $E$ and $\Theta : \R^{E} \to \R^{|E|}$ denotes the usual valuation map
$$\Theta(g) = (g(e_1),\ldots,g(e_{|E|})),$$
we can naturally identify $P$ with a subset of points in $\mathbb{R}^{|E|}$, its image $\Theta(P)$. Since $G$ is bipartite, this set of points has a remarkable characterization. 
For each matching $M \in \mathcal{M}_{m}$, let us denote by $\chi_{M} \in \mathbb{R}^{E}$ the characteristic function of $M$ as a subset of $E$, namely
$$\chi_{M}(e)=
	\begin{cases}
		1 \text{ if } e\in M \\
		0\text{ otherwise.} 
	\end{cases}$$
Furthermore, let $\overrightarrow{\chi_{M}}$ denote the point $(\chi_{M}(e))_{e \in E} \in \mathbb{R}^{|E|}$. Using the standard ``cycle and path elimination" proof of \cite[Theorem 7.1.2, page 269]{Lovasz}, it is not difficult to see that $\Theta(P)$ must be precisely the convex hull of $\left\{\overrightarrow{\chi_{M}}: M \in \mathcal{M}_{m}\right\}$. Indeed, first of all, note that if $\alpha_{M} \in [0,1]$ are such that $\sum_{M \in \mathcal{M}_{m}} \alpha_{M} = 1$, then for every $v \in A$ we have
$$\sum_{u \in B} \sum_{M \in \mathcal{M}_{m}} \alpha_{M}\chi_{M}(v,u) =  \sum_{M \in \mathcal{M}_{m}} \alpha_{M} \sum_{v \in A} \chi_{M}(v,u) \leq \sum_{M \in \mathcal{M}_{m}} \alpha_{M} = 1.$$
Similarly, for every $u \in B$, we have $\sum_{v \in A} \sum_{M \in \mathcal{M}_{m}} \alpha_{M}\chi_{M}(v,u) \leq 1$, and also
$$\sum_{v \sim u} \sum_{M \in \mathcal{M}_{m}} \alpha_{M}\chi_{M}(v,u) = \sum_{M \in \mathcal{M}_{m}} \alpha_{M} \sum_{v \sim u} \chi_{M}(v,u) = m \sum_{M \in \mathcal{M}_{m}} \alpha_{M} = m.$$
This shows that the convex hull of $\left\{\overrightarrow{\chi_{M}}: M \in \mathcal{M}_{m}\right\}$ is contained in $\Theta(P)$. To show the reverse, the idea is to consider each extremal point in $\Theta(P)$ and show that it must be the characteristic vector of a matching of size $M$. In other words, if $g \in P$ is such that $\Theta(g)$ is an extremal point in $\Theta(P)$, then it suffices to show that the support of $g$ must always be a matching of size $m$. To see this, one proceeds by contradiction: if the support of $g$ is not a matching of size $m$, then one can always produce (explicitly) some $g_1,\ldots,g_k \in P$ and $c_1,\ldots,c_k \in [0,1]$ with $\sum_{i=1}^{k} c_i = 1$ and such that $$g = \sum_{i=1}^{k} c_i g_i.$$
These auxiliary fractional matchings can be found algorithmically by analyzing the structure of the graph on the underlying set of edges in the support of $g$. We skip the details here, since the analysis is quite similar to the argument for the (full) matching polytope. For more details, see \cite[Chapter 18]{Schrijver}.

Given any function $f \in P$, we now know that there always exist coefficients $\mu_{M} \in [0,1]$, $M \in \mathcal{M}_{m}$, such that $\sum_{M \in M_{m}} \mu_{M} = 1$ and 
$$\Theta(f) = \sum_{M \in M_{m}} \mu_{M} \overrightarrow{\chi_{M}}.$$
Pulling this back to $\mathbb{R}^{E}$, it is easy to see that the coefficients of this affine combination give us the desired probability distribution.  
\end{proof}

\paragraph{Weighted comparability graphs.} Suppose that $i \le \frac{n-1}{2} D$. Let us consider a weighted bipartite graph $G_i$ between $V_i$ and $V_{i+1}$ in which two vectors $v \in V_i$ and $u \in V_{i+1}$ are connected by an edge if $v \prec u$. For any edge $(v, u)$ there exists a unique coordinate $t \in [D]$ such that $u_t = v_t+1$. In this case, we assign color $u_t$ to the edge $(v, u)$ and define $w_{vu} = u_t(n-u_t)$.
For each $v \in V_i$ above, its (weighted) degree in $G_i$ is given by 
$$
\deg_{G_i} v := \sum_{(v,u) \in E} w_{v u}.$$
Similarly, for each $u \in V_{i+1}$ we have 
$$
\deg_{G_i} u := \sum_{(v,u) \in E} w_{v u}.$$
Finally, for every $i < \frac{n-1}{2} D$, we say that a matching between $V_{i}$ and $V_{i+1}$ is {\it{full}} if its edges cover all the vertices from the set $V_i$. 
We prove the following important preliminary result. 

\begin{lemma} \label{distributions}
For $i \le \frac{n-1}{2}D - \frac{n-1}{2}$ there exists a distribution $\mu_i$ on the set of full matchings between $V_i$ and $V_{i+1}$ such that each edge $(v, u) \in V_{i} \times V_{i+1}$ belongs to a $\mu_i$-random matching with probability $\frac{w_{v u}}{\deg_{G_i} v}$.

For $\frac{n-1}{2}D - \frac{n-1}{2} < i \le \frac{n-1}{2} D$ there exists a distribution $\mu_i$ on the set of matchings of size at least $|V_i| (1 - 3/D)$ between $V_i$ and $V_{i+1}$ such that each edge $(v, u) \in V_{i} \times V_{i+1}$ belongs to a $\mu_i$-random matching with probability at most $\frac{w_{v u}}{\deg_{G_i} v}$.

For $i > \frac{n-1}{2} D$ the analogous results hold, since we can define the corresponding distribution $\mu_i$ by symmetry.
\end{lemma}

\begin{proof} Since the weights are explicit, we can actually compute these degrees as follows. 
For $v \in [n]^D$ and $j = 0, \ldots, n-1$ let $b_j(v)$ be the number of coordinates of $v$ equal to $j$. Note that for $v \in V_i$ we have
\begin{align}\label{123}
    \sum_{j = 0}^{n-1} b_j(v) = D,~~    \sum_{j = 0}^{n-1} b_j(v) j = i.
\end{align}
For convenience, denote $a_j = j(n-j)$ for each $j=0,\ldots,n-1$. Note that $a_{j+1} = a_j + n- 2j-1$ and that $a_0 = a_{n} = 0$.

Then for each $v \in V_{i}$, we have
\begin{align}
\deg_{G_i} v = \sum_{(v,u) \in E} w_{v u} = \sum_{j = 0}^{n-2} b_j(v) a_{j+1},
\end{align}

For $u \in V_{i+1}$, we analogously have
\begin{align}
\deg_{G_i} u = \sum_{(v,u) \in E} w_{v u} = \sum_{j = 1}^{n-1} b_j(u) a_{j}.
\end{align}

Furthermore, note that if $(v, u)$ is an edge of color $l$ (i.e. $u_t = v_t+1 = l$ for some $t\in [D]$) then $b_l(v) = b_l(u) - 1, b_{l-1}(v) = b_{l-1}(u)+1$ and $b_{j}(v) = b_j(u)$ for $j \neq l, l-1$ and so
\begin{eqnarray*}
\deg_{G_i} (v) &=& a_l -a_{l+1} + \sum_{j = 0}^{n-2} b_j(u) a_{j+1} \\
&=& 2l -n+1 + b_0(u) a_1 + \sum_{j = 1}^{n-1} b_j(u) a_{j+1} \\
&\stackrel{(\ref{123})}{=}& 2l - n + 1 + b_0(u) a_1 + \deg_{G_i}(u) + (D - b_0(u)) (n-1) - 2(i+1) \\
&=&\deg_{G_i}(u) + D (n-1) - 2(i+1) + 2 l - n+1.
\end{eqnarray*}

Denote $\delta = D (n-1) - 2(i + 1) \ge 0$.
Consider the function $f(v, u) = \frac{w_{v u}}{\deg_{G_i} v}$. Then for fixed $u \in V_{i+1}$ we need to upper bound the following:
\begin{align}\label{f}
    \sum_{(v,u) \in E} \frac{w_{v u}}{\deg_{G_i} v} = \sum_{j = 0}^{n-1}  \frac{b_j(u) a_j}{\deg_{G_i}(u) + \delta + 2j - n+1}.
\end{align}

Note that if $\delta \ge n-1$ then $\delta + 2j - n+1 \ge 0$ and so the sum can be bounded from above by
$$
\sum_{j = 0}^{n-1} \frac{b_j(u) a_j}{\deg_{G_i}(u)} = 1.
$$
This means that for $\delta \ge n-1$ the function $f$ satisfies the conditions of Lemma \ref{polytope}, and so the first conclusion of Lemma \ref{distributions} follows.

Now we suppose that $0 \le \delta < n-1$. It is then clear that for any $u \in V_{i+1}$ has at least $\frac{D-1}{2}$ non-zero coordinates and we have 
\begin{equation}\label{deg}
\deg_{G_i}u = \sum_{j = 0}^{n-1} b_j(u) a_j \ge \frac{(D-1)(n-1)}{2}.
\end{equation}
So for any $v \preceq u$, $v \in V_i$ we have
$$
\frac{w_{v u}}{\deg_{G_i} v} \le \frac{w_{v u}}{\deg_{G_i} u - n+1} \le \frac{D-1}{D - 3} \frac{w_{v u}}{\deg_{G_i} u}.
$$
We conclude that 
\begin{equation}\label{asda}
\sum_{(v,u) \in E} \frac{w_{v u}}{\deg_{G_i} v} \le \frac{D-1}{D-3}.
\end{equation}
Let $\theta \le \frac{D-3}{D-1}$ be the largest number such that $\theta |V_{i}|$ is an integer. It is easy to check that $\theta \ge 1 - 3/D$. Then (\ref{asda}) implies that the function $g(v, u) = \theta f(v, u)$ satisfies the conditions of Lemma \ref{polytope} with $m = \theta |V_{i}| \ge (1-3/D) |V_i|$ so the second assertion of Lemma \ref{distributions} follows.


\end{proof}

We are now ready to prove Theorem \ref{super}. 

\begin{proof}
For each $i = 0, \ldots, (n-1)D-1$, consider the probability distribution $\mu_i$ from Lemma \ref{distributions}, and let $M_i$ be a $\mu_i$-random matching between $V_i$ and $V_{i+1}$ (possibly partial, depending on how close we are to the middle layer).
The union of these matchings defines a graph on $[n]^D$ whose connected components are chains with respect to the partial order on $[n]^D$. 
Let $\{C_1, \ldots, C_N\}$ be the set of these chains. It is perhaps worth emphasizing that in this list some of these chains may be trivial, in the sense that a $C_i$ could in principle consist of one single element of $[n]^{D}$ (if this element was not covered by some matching). Consequently, we have a chain decomposition $[n]^D = \bigsqcup_{j = 1}^N C_j$. Nevertheless, we claim that this decomposition satisfies
\begin{equation} \label{chains}
N \le \left(1 + \frac{3n}{D}\right) N_{n, D}.
\end{equation}
To prove this estimate, we can argue as follows. Denote $k = \lfloor \frac{n-1}{2} D\rfloor$. 
For each set $C_j$ let $v_j \in C_j$ be the closest element to the middle layer $V_k$. The fact that for $i \not \in [k - \frac{n-1}{2}, k + \frac{n+1}{2}]$ the matching $M_i$ completely covers the smaller part of the bipartition $(V_i, V_{i+1})$ implies that $v_j \in V_{x_j}$ for some 
$x_j \not \in [k - \frac{n-1}{2}, k + \frac{n+1}{2}]$ and all $j$. On the other hand, the fact that $|M_i| \ge (1-3/D)\min\{|V_i|, |V_{i+1}|\} $ implies that for any $i \neq k$ the number of $j$ for which $x_j = i$ is at most $\frac{3}{D} \min\{|V_i|, |V_{i+1}|\} \le \frac{3}{D} N_{n, D}$. The bound on $N$ follows immediately. 

Now, suppose that $X \subset [n]^D$ is a set of size at least $a(1 + 3n/D) N_{n, D} + b$. By \eqref{chains} and the pigeonhole principle, we know that there are always at least $b$ pairs $(x, y)$ of elements of $X$ belonging to the same set $C_j$, such that $x \prec y$ and so that the rank of $y$ exceeds the rank of $x$ by at least $a$. Denote by $\xi$ the random variable equal to the total number of such pairs $(x, y)$. 

On the other hand, let us fix such a pair $x \prec y$ with $x, y \in X$ and estimate the probability that $x, y$ belong to the same set $C_j$ for some $j$. 
Assume that $x \in V_i$ and $y \in V_{i+a}$. Let $i+j \in [i, i+a]$ be the closest index to $k$. A chain between $x$ and $y$ has the form 
$$x = z_0 \prec z_1 \prec \ldots \prec z_a = y,$$
and there are at most $a!$ ways to fix such a chain. By the independence of matchings $M_l$ and the choice of $\mu_l$, the probability that $\{z_0, \ldots, z_a\}$ is contained in a set $C_j$ for some $j$ is at most $f(z_0, z_1) f(z_1, z_2) \ldots f(z_{a-1}, z_a)$. Fix $0\le j' < j$, by a computation analogous to (\ref{deg}), we have 
$$
f(z_{j'}, z_{j'+1}) \le \frac{w_{z_{j'}, z_{j'+1}}}{\deg_{G_{i+j'}}z_{j'}} \le \frac{n^2 /4}{(n-1)D / 2} \le \frac{n}{D},
$$
and similarly for $j' \ge j$ one can show that $f(z_{j'}, z_{j'+1}) \le \frac{n}{D}$. So the probability that $z_0, \ldots, z_a$ are contained in the same chain is at most $(n/D)^a$. Since the number of chains connecting $x$ and $y$ is at most $a!$, the desired probability is at most $a! (n/D)^a$.



We conclude that 
$$
\E \xi \le a! \left(\frac{n}{D}\right)^a \# \{x \preccurlyeq y \in X,~ r(y) - r(x) \ge a\}.
$$
But $\xi \ge b$ always so 
$$
\# \{x \preccurlyeq y \in X,~ r(y) - r(x) \ge a\} \ge \frac{b}{a! n^a}  D^a.
$$
This concludes the proof of Theorem \ref{super}. 
\end{proof}

\subsection{Application of container lemma} 

With the result from Theorem \ref{super}, we can now apply the graph container method in order to estimate the number of antichains in $[n]^{D}$. More precisely, we shall take advantage of \cite[Lemma 5.4]{NSS}, which restate below for convenience. 

\begin{lemma} \label{containers}
For $k \geq 1$, let $d_1 > \ldots > d_k$ and $m_0 > m_1 > \ldots > m_k$ be positive integers and let $P$ be a poset such that $|P|=m_0$, and, for $1 \leq j \leq k$, every subset $S$ of $P$ of cardinality greater than $m_j$ contains at least $|S|d_j$ comparable pairs. Then, there is a collection $\mathcal{F}$ of subsets of $P$ such that
\begin{enumerate}
\item $|\mathcal{F}| \leq \prod_{r=1}^{k}{m_{r-1} \choose \leq m_{r-1}/(2d_{r}+1)}$,
\item $|A| \leq m_{k}+\sum_{r=1}^{k}\frac{m_{r}-1}{2d_{r}+1}$ holds for every $A \in \mathcal{F}$, and
\item for every antichain $I$ of $P$, there exists $A \in \mathcal{F}$ such that $I \subset A$. 
\end{enumerate}
\end{lemma}

We apply Lemma \ref{containers} when $P$ is our poset $[n]^{D}$ and for $k=2$. In this case, we clearly have $m_{0} = n^{D}$. We then let $d_1 := \frac{D}{2n}$, $d_2 := \frac{\sqrt{D}}{2n}$, $m_1 := N_{n, D}(2 + 6n/D)$, and $m_2 := N_{n, D} (1 + 3n/D +1/\sqrt{D})$. With this choice of parameters, Theorem \ref{super} implies that every set $X \subset [n]^D$ of size at least $m_j$ contains at least $d_j |X|$ comparable pairs, so Lemma \ref{containers} applies. We thus get a family of containers $\mathcal{F}$ for the set of antichains in $[n]^D$ such that the size of $\mathcal{F}$ satisfies
$$|\mathcal{F}| \le {n^D \choose \le \frac{n^D}{2d_1 + 1}} {3 N_{n, D} \choose \le \frac{3 N_{n, D}}{2 d_2 + 1}} \le (6 d_1)^{n^D / d_1} (6 d_2)^{N_{n, D} / d_2} \le 2^{c_n N_{n, D} \log D / \sqrt{D}}.$$
Here $c_{n}$ is a positive constant which depends solely on $n$. In the first inequality, we used the general fact that ${m \choose \le t} \leq (em/t)^t$ holds for all $1 \leq t \leq m$, whereas in the last inequality we used the fact that $N_{n,D} = O(n^{D}/\sqrt{D})$, which follows immediately from (the renormalized version of) the estimate from \eqref{CLT}. 

Moreover, since $m_1 = N_{n, D}(2 + 6n/D) \le 3 N_{n, D}$ and 
$$m_2 = N_{n, D} \left(1 + \frac{3n}{D} + \frac{1}{\sqrt{D}}\right) \le N_{n, D} \left(1 + \frac{2}{\sqrt{D}}\right),$$
every container $A \in \mathcal{F}$ satisfies
$$|A| \leq m_2 + \frac{m_1}{2d_2 + 1} \le N_{n, D} \left(1 + \frac{2}{\sqrt{D}}\right) + \frac{3 N_{n, D}}{2d_2+1} \leq N_{n, D} \left(1 + \frac{c'_n}{\sqrt{D}}\right),$$
for some positive constant $c'_n$ which also depends only on $n$. Last but not least, by design, the collection of subsets $\mathcal{F}$ is such that for any antichain $X \subset [n]^D$ there exists a set $A \in \mathcal{F}$ such that $X \subset A$.

By a simple union bound, we can now conclude that the total number $\mathcal{A}_{D,n}$ of antichains in $[n]^D$ satisfies
$$\mathcal{A}_{D,n} \leq 2^{N_{n, D} (1 + c'_n / \sqrt{D}) + c_n N_{n, D} \log D / \sqrt{D}} \le 2^{\left(1 + \frac{C_n \log D}{\sqrt{D}}\right) N_{n, D}},
$$
for some positive (new) constant $C_n$ that depends solely on $n$. This completes the proof of \eqref{anti2} and thus that of Theorem \ref{anti}.

\bigskip

{\bf{Acknowledgements}}. We would like to thank J\'ozsef Balogh, Richard Kenyon, and Igor Pak for helpful discussions. We would also like to thank the anonymous referees for useful comments which improved the final version of the paper.

\end{document}